\documentclass[a4paper,12pt]{article}

\usepackage[english]{babel}
\usepackage{amsfonts,amsmath,amsthm,mathrsfs,amscd,amssymb}
\usepackage{cite}
\usepackage{authblk}
\usepackage{mhchem}

\newtheorem{T}{Theorem}

\newtheorem{prop}{Proposition}

\theoremstyle{remark}

\providecommand{\keywords}[1]
{
  \small	
  \textit{Key Words:~} #1
}

\providecommand{\MSC}[1]
{
  \small	
  \textit{Mathematics Subject Classification 2010:~} #1
}

\begin{document}

\title{On the characterization of a finite random field by conditional distribution and its Gibbs form}

\author[1]{Khachatryan L.A.}
\author[2]{Nahapetian B.S.}

\affil[1]{\small Institute of Mathematics, NAS of RA, {\it linda@instmath.sci.am}}
\affil[2]{\small Institute of Mathematics, NAS of RA, {\it nahapet@instmath.sci.am}}

\date{}

\maketitle

\begin{abstract}
In this paper, we show that the methods of mathematical statistical physics can be successfully applied to random fields in finite volumes. As a result, we obtain simple necessary and sufficient conditions for the existence and uniqueness of a finite random field with a given system of one-point conditional distributions. Using the axiomatic (without the notion of potential) definition of Hamiltonian, we show that any finite random field is Gibbsian. We also apply the proposed approach to Markov random fields.
\end{abstract}


\MSC{60G60, 60E05, 60J99, 62H99, 82B03}

\section{Introduction}

The idea of characterization of a multivariate probability distribution by means of conditional probabilities is very old and goes back to the concept of the finite Markov chain. Already here, from the very beginning, we are faced with two fundamentally important characterization problems: the problem of restoring a chain from the system of its transition probabilities and the problem of the existence of a Markov chain corresponding to a given stochastic matrix.

Both of these problems are relevant in the general setting, when one deals with an arbitrary random field in an infinite or finite volume. Having solved the problem of restoring a random field from the system of its conditional probabilities, we can establish the general properties of the random field in these terms. The existence of a random field with a given system of conditional distributions allows constructing models of random fields with the required properties.

In mathematical statistical physics, as a rule, probabilistic models in infinite volumes are studied. This is due to the fact that it is in such models that the phenomenon of phase transition is best described (see, for example, the book of Sinai~\cite{Sinai}).

Dobrushin~\cite{Dobr} was the first to consider the problem of characterization of a random field on an infinite integer lattice by specification --- a system of consistent probability distributions parameterized by infinite boundary conditions. He obtained criteria for both the existence and uniqueness of a random field with a given specification. The basic example of a specification is the Gibbs specification. Further, Dobrushin~\cite{Dobr-Gibbs} introduced the concept of a Gibbs random field as a field the conditional distribution of which almost everywhere coincides with the Gibbs specification. This concept was later independently discovered by Lanford and Ruelle~\cite{LanfRu}.

There is a large number of works devoted to the characterization problems of a random field specified in a finite volume (finite random field). The main issue considered here is to find the conditions under which a given system of probability distributions is compatible, that is, there exists a random field which conditional probabilities coincide with the elements of the given system. This issue was studied by Besag~\cite{Besag}, Abrahams and Thomas~\cite{AbTh}, Arnold and Press~\cite{ArPr}, Gurevich~\cite{Gurev}, and Chen~\cite{Chen} among others. For finite random fields, it is also of interest to find the conditions under which a given field is Gibbsian. This question was considered, for example, by Hammersley and Clifford~\cite{HammCliff}, Griffeath~\cite{KSK} and Kaiser and Cressie~\cite{KC}.

The methods used to solve the characterization problems in finite and infinite volumes are fundamentally different. Apparently, this led to the situation where the research in these two directions was carried out independently. At the same time, the methods developed in these areas can be mutually beneficial. As an example, we cite the Grimmett's theorem~\cite{Grimmet} on the construction of the interaction potential of a finite Gibbs random field based on the M\"{o}bius transform. This approach was used by Sullivan~\cite{Sull} and Kozlov~\cite{Kozlov} to find the conditions under which the Dobrushin's specification has a Gibbs form with some potential. Representation of a specification or a finite random field in Gibbs form allows one to study their properties using the well-developed mathematical apparatus of the theory of Gibbs random fields. This approach is especially productive when studying Markov (see, for example, works by Griffeath~\cite{KSK}, Sullivan~\cite{Sull}, Spitzer~\cite{Spitzer}, and Georgii~\cite{Georgii}) and Gaussian random fields (we refer Dobrushin~\cite{DGauss}, Georgii~\cite{Georgii} and K\"{u}nsch~\cite{Kunsch} among others).

In this paper, we show that the methods developed in~\cite{DN04} and~\cite{DN19} for infinite random fields can be successfully applied to random fields in finite volumes. As a result, we obtain simple necessary and sufficient conditions for the existence and uniqueness of a finite random field with a given system of one-point conditional distributions. In this case, we establish an explicit formula for the multivariate distribution in terms of the given one-point conditional distributions. In addition, it is shown that in a finite volume the Gibbs form is universal, namely, the elements of any compatible system of probability distributions have the Gibbs form. As an application, we obtain the generalization of the well-known Hammersley--Clifford theorem~\cite{HammCliff} on the equivalence of Markov and Gibbs finite random fields.

For simplicity, in this article, we consider random fields with a finite state space. However, all the results can be extended to the case of both countable and continuous state spaces.

The paper has the following structure. In Section 2, we provide basic notations. Section 3 is devoted to the problem of characterization of finite random fields by systems of consistent one-point distributions parameterized by boundary conditions. In Section 4, we consider the question of the Gibbsianess of finite random fields. Section 5 provides some concluding remarks.

\section{Notations}

Let $\Lambda$ be a set with a finite number of elements, $\vert \Lambda \vert < \infty$, and let each point $t \in \Lambda$ be associated with the set $X^t$, which is a copy of some finite set $X$. For any $V \subseteq \Lambda$, denote by $X^V = \{  x = (x_t, t \in V): \, x_t \in X^t, t \in V\}$ the set of functions (configurations) defined on $V$ and tacking values in $X$.

For any $V \subset \Lambda$, denote by $x_V$ the restriction of configuration $x \in X^\Lambda$ on $V$. When denoting one-point sets $\{t\}$, $t \in \Lambda$, braces will be omitted in some cases. In particular, for the restriction of a configuration $x \in X^\Lambda$ on $t \in \Lambda$, the notation $x_t$ will be used instead of $x_{\{t\}}$. For any $V, I \subset \Lambda$ such that $V \cap I = \emptyset$ and any $x \in X^V$, $y \in X^I$, denote by $xy$ the concatenation of $x$ with $y$, that is, the configuration on $V \cup I$ equal to $x$ on $V$ and to $y$ on $I$. The concatenation of $x$ with an empty configuration $\boldsymbol{\emptyset}$ we assume to coincide with $x$, i.e., $x \boldsymbol{\emptyset} = x$ for all $x \in X^V$, $V \subseteq \Lambda$.

Probability distribution $P_\Lambda$ on $X^\Lambda$ will be called a (\emph{finite}) \emph{random field}. A random field $P_\Lambda$ is called \emph{positive} if $P_\Lambda(x)>0$ for all $x \in X^\Lambda$. By $P_V$ we denote the restriction of $P_\Lambda$ on $V \subset \Lambda$, that is,
$$
P_V(x) = \sum \limits_{u \in X^{\Lambda \backslash V}} P_\Lambda(xu), \qquad x \in X^V.
$$
For $V = {\emptyset}$, one has $P_{\emptyset}(\boldsymbol{\emptyset}) = 1$.

For a random field $P_\Lambda$, define a conditional probability distribution
$$
Q_V^z(x) = \frac{P_{V \cup I}(xz)}{P_I(z)}, \qquad x \in X^V, z \in X^I,
$$
on $X^V$ under the (boundary) condition $z$, where $V,I \subset \Lambda$, $V \cap I = \emptyset$ and $P_I(z) > 0$.

\section{Characterization of random fields by conditional distributions}
\label{Sec-Q->P}

In this section, we give necessary and sufficient conditions under which there is a unique random field, which conditional probabilities coincide with the corresponding elements of a given system of probability distributions.

\subsection{Positive random fields}

Let $P_\Lambda$ be a positive random field and
\begin{equation}
\label{P->Q1}
Q_t^z(x) = \frac{P_\Lambda(xz)}{\sum \limits_{u \in X^t} P_\Lambda(uz)}, \qquad x \in X^t, z \in X^{\Lambda \backslash \{t\}}, t \in \Lambda,
\end{equation}
be its one-point conditional probabilities. The family $\boldsymbol{Q}_1(P_\Lambda) = \{Q_t^z, z \in X^{\Lambda \backslash \{t\}}, t \in \Lambda \}$ we will call \emph{one-point conditional distribution} of the random field $P_\Lambda$.

Let $\boldsymbol{\rm Q}_1 = \{q_t^z, z \in X^{\Lambda \backslash \{t\}}, t \in \Lambda \}$ be a family of one-point probability distributions $q_t^z$ on $X^t$ parameterized by configurations $z \in X^{\Lambda \backslash \{t\}}$, $t \in \Lambda$. In mathematical statistical physics, the parameter $ z $ is usually called the boundary condition. In the future, we will adhere to this terminology.

A random field $P_\Lambda$ is said to be \emph{compatible} with a system $\boldsymbol{\rm Q}_1$ if $\boldsymbol{Q}_1 (P_\Lambda) = \boldsymbol{\rm Q}_1$.

\begin{T}
\label{Th-Q1->P+}
For a given system $\boldsymbol{\rm Q}_1 = \{q_t^z, z \in X^{\Lambda \backslash \{t\}}, t \in \Lambda \}$ of positive one-point probability distributions parameterized by boundary conditions, there exists a unique random field $P_\Lambda$ compatible with it if and only if the elements of $\boldsymbol{\rm Q}_1$ satisfy the following consistency conditions: for all $t,s \in \Lambda$ and $x,u \in X^t$, $y,v \in X^s$, $z \in X^{\Lambda \backslash \{t,s\}}$, it holds
\begin{equation}
\label{Q1-sogl}
q_t^{zy} (x) q_s^{zx} (v) q_t^{zv} (u) q_s^{zu} (y) = q_t^{zy} (u) q_s^{zu} (v) q_t^{zv} (x) q_s^{zx} (y).
\end{equation}
\end{T}
\begin{proof}
The necessity of conditions~\eqref{Q1-sogl} follows directly from the definition of conditional probabilities. Let us prove the sufficiency by giving an explicit form of the finite random field $P_\Lambda$ defined by the system $\boldsymbol{\rm Q}_1$.

For any $x \in X^\Lambda$, put
\begin{equation}
\label{Q1->P}
P_\Lambda(x) = \prod \limits_{j=1}^n \frac{q_{t_j}^{(xu)_j} (x_j)}{q_{t_j}^{(xu)_j} (u_j)} \cdot \left( \sum \limits_{\alpha \in X^\Lambda} \prod \limits_{j=1}^n \frac{q_{t_j}^{(\alpha u)_j} (\alpha_j)}{q_{t_j}^{(\alpha u)_j} (u_j)} \right)^{-1},
\end{equation}
where $u \in X^\Lambda$, $\Lambda = \{t_1, t_2, ..., t_n\}$ is some enumeration of points in $\Lambda$, $n = \vert \Lambda \vert$, and
$$
(xu)_j = x_1 ... x_{j-1} u_{j+1} ... u_n, \quad 1 < j < n, \qquad (xu)_1 = u_2 u_3 ... u_n, \quad (xu)_n = x_1 x_2 ... x_{n-1},
$$
$x_j = x_{t_j}$, $u_j = u_{t_j}$, $1 \le j \le n$. Expanded, we have
$$
\prod \limits_{j=1}^n \frac{q_{t_j}^{(xu)_j} (x_j)}{q_{t_j}^{(xu)_j} (u_j)} = \dfrac{q_{t_1}^{u_2 u_3 ... u_n}(x_1)}{q_{t_1}^{u_2 u_3 ... u_n}(u_1)} \cdot \frac{q_{t_2}^{x_1 u_3 ... u_n}(x_2)}{q_{t_2}^{x_1 u_3 ... u_n}(u_2)} \cdot \frac{q_{t_3}^{x_1 x_2 u_4 ... u_n}(x_3)}{q_{t_3}^{x_1 x_2 u_4 ... u_n}(u_3)} \cdot ... \cdot \frac{q_{t_n}^{x_1 x_2 ... x_{n-1}}(x_n)}{q_{t_n}^{x_1 x_2 ... x_{n-1}}(u_n)}.
$$
First of all, let us make sure that the formula~\eqref{Q1->P} is correct. To do this, it is necessary to check that the values of $P_\Lambda$ do not depend on the enumeration of the elements of $\Lambda$ and on the choice of the configuration $u \in X^\Lambda$.

We start by showing that for fixed $x,u \in X^\Lambda$, the right-hand side of~\eqref{Q1->P} does not depend on the enumeration of the elements of $\Lambda$. Note that any enumeration of the elements of $\Lambda$ can be obtained from another one by successive changing the numeration of two points. For some $k$, $1 < k \le n$, consider the numeration $\varphi(j)$, $1 \le j \le n$, such that $t_{\varphi(k-1)} = t_k$, $t_{\varphi(k)} = t_{k-1}$, and $t_{\varphi(j)} = t_j$ for $j \neq k -1$ and $j \neq k$. Then we have
\begin{eqnarray*}
\prod \limits_{j=1}^n \dfrac{q_{t_{\varphi(j)}}^{(xu)_{\varphi(j)}} (x_{\varphi(j)})}{q_{t_{\varphi(j)}}^{(xu)_{\varphi(j)}} (u_{\varphi(j)})} & = & \prod \limits_{j\neq k-1,k} \frac{q_{t_j}^{(xu)_j} (x_j)}{q_{t_j}^{(xu)_j} (u_j)} \\
\\
&& \cdot \dfrac{ q_{t_k}^{x_1 ... x_{k-2} u_{k-1} u_{k+1}... u_n} (x_k)}{ q_{t_k}^{x_1 ... x_{k-2} u_{k-1} u_{k+1} ... u_n} (u_k)} \cdot \dfrac{q_{t_{k-1}}^{x_1 ... x_{k-2} x_k u_{k+1} ... u_n} (x_{k-1}) }{q_{t_{k-1}}^{x_1 ... x_{k-2} x_k u_{k+1} ... u_n} (u_{k-1}) }.
\end{eqnarray*}
Denoting  $z = x_1 ... x_{k-2} u_{k+1}... u_n$, we have
$$
\dfrac{ q_{t_k}^{x_1 ... x_{k-2} u_{k-1} u_{k+1}... u_n} (x_k)}{ q_{t_k}^{x_1 ... x_{k-2} u_{k-1} u_{k+1} ... u_n} (u_k)} \cdot \dfrac{q_{t_{k-1}}^{x_1 ... x_{k-2} x_k u_{k+1} ... u_n} (x_{k-1}) }{q_{t_{k-1}}^{x_1 ... x_{k-2} x_k u_{k+1} ... u_n} (u_{k-1}) } = \dfrac{ q_{t_k}^{z u_{k-1}} (x_k)}{ q_{t_k}^{z u_{k-1}} (u_k) } \cdot \dfrac{q_{t_{k-1}}^{z x_k} (x_{k-1}) }{q_{t_{k-1}}^{z x_k} (u_{k-1}) }.
$$
According to the consistency conditions~\eqref{Q1-sogl} of the elements of the system $\boldsymbol{\rm Q}_1$, we can write
$$
\dfrac{ q_{t_k}^{z u_{k-1}} (x_k)}{ q_{t_k}^{z u_{k-1}} (u_k) } \cdot \dfrac{q_{t_{k-1}}^{z x_k} (x_{k-1}) }{q_{t_{k-1}}^{z x_k} (u_{k-1}) } = \dfrac{ q_{t_{k-1}}^{z u_k} (x_{k-1})}{q_{t_{k-1}}^{z u_k} (u_{k-1}) } \cdot \dfrac{q_{t_k}^{z x_{k-1}} (x_k) }{ q_{t_k}^{z x_{k-1}} (u_k)}.
$$
Noting that
$$
\dfrac{ q_{t_{k-1}}^{z u_k} (x_{k-1})}{q_{t_{k-1}}^{z u_k} (u_{k-1}) } \cdot \dfrac{q_{t_k}^{z x_{k-1}} (x_k) }{ q_{t_k}^{z x_{k-1}} (u_k)} = \dfrac{q_{t_{k-1}}^{(xu)_{k-1}}(x_{k-1})}{q_{t_{k-1}}^{(xu)_{k-1}}(u_{k-1}) } \cdot \dfrac{q_{t_k}^{(xu)_k}(x_k)}{q_{t_k}^{(xu)_k}(u_k)},
$$
we conclude
$$
\prod \limits_{j=1}^n \dfrac{q_{t_{\varphi(j)}}^{(xu)_{\varphi(j)}} (x_{\varphi(j)})}{q_{t_{\varphi(j)}}^{(xu)_{\varphi(j)}} (u_{\varphi(j)})} = \prod \limits_{j=1}^n \frac{q_{t_j}^{(xu)_j} (x_j)}{q_{t_j}^{(xu)_j} (u_j)}.
$$
%
%or, equivalently,
%$$
%\frac{ q_{t_{k-1}}^{z u_k} (x_{k-1}) q_{t_k}^{z x_{k-1}} (x_k)}
%  { q_{t_{k-1}}^{z u_k} (u_{k-1}) q_{t_k}^{z x_{k-1}} (u_k) }
%=
%\frac{ q_{t_k}^{z u_{k-1}} (x_k) q_{t_{k-1}}^{z x_k} (x_{k-1})}
%  { q_{t_k}^{z u_{k-1}} (u_k) q_{t_{k-1}}^{z x_k} (u_{k-1})},
%$$
%where $z = x_1 ... x_{k-2} u_{k+1} ... u_n$, the value of the products~\eqref{product-fi} and~\eqref{product} are equal.

In order to prove that $P_\Lambda$ does not depend on the choice of the configuration $u \in X^\Lambda$, it is sufficient to verify that the right-hand side of~\eqref{Q1->P} will not change when substituting some other configuration $y \in X^\Lambda$ instead of $u$. Note that~\eqref{Q1->P} can be written in the following form
$$
P_\Lambda(x) = q_{t_1}^{u_{\Lambda \backslash \{t_1\}}} (x_1) \prod \limits_{j=2}^n \frac{q_{t_j}^{(xu)_j} (x_j)}{q_{t_j}^{(xu)_j} (u_j)} \left( \sum \limits_{\alpha \in X^\Lambda} q_{t_1}^{u_{\Lambda \backslash \{t_1\}}} (\alpha_1) \prod \limits_{j=2}^n \frac{q_{t_j}^{(\alpha u)_j} (\alpha_j)}{q_{t_j}^{(\alpha u)_j} (u_j)} \right)^{-1},
$$
and hence, the value of $P_\Lambda$ does not depend on the $t_1$ component of the configuration $u$. Since it also does not depend on the enumeration of the elements of $\Lambda$, the value of $P_\Lambda$ does not change if we substitute into~\eqref{Q1->P} a configuration $u' \in X^\Lambda$ that differs from $u$ at one component. It remains to note that any configuration $y$ can be obtained from configuration $u$ by sequentially changing one component.

It is not difficult to see that $P_\Lambda$ is a probability distribution on $X^\Lambda$, and $P_\Lambda(x) > 0$ for all $x \in X^\Lambda$. To conclude the proof of the theorem, it remains to verify that $\boldsymbol{Q}_1(P_\Lambda) = \boldsymbol{\rm Q}_1$. For this, we need to show that $Q_t^z(x) = q_t^z(x)$ for any $t \in \Lambda$ and $x \in X^t$, $z \in X^{\Lambda \backslash \{t\}}$.

Since the values of $P_\Lambda $ do not depend on the enumeration of the elements in $\Lambda$, for every $t \in \Lambda$ and all $x \in X^t$, $z \in X^{\Lambda \backslash \{t \}}$, we can assume $t = t_n$ and write
\begin{eqnarray*}
Q_t^z(x) & = & \dfrac{P(xz)}{\sum \limits_{\alpha \in X^t} P(\alpha z)} = \dfrac{q_{t_1}^{(zu)_1}(z_1) q_{t_2}^{(zu)_2}(z_2) \cdot ... \cdot q_{t_{n-1}}^{(zu)_{n-1}}(z_{n-1}) q_{t_n}^{(zu)_n}(x)}{q_{t_1}^{(zu)_1}(u_1) q_{t_2}^{(zu)_2}(u_2) \cdot ... \cdot q_{t_{n-1}}^{(zu)_{n-1}}(u_{n-1}) q_{t_n}^{(zu)_n}(u_n)} \\
& & \cdot \left( \sum \limits_{\alpha \in X^t} \dfrac{q_{t_1}^{(zu)_1}(z_1) q_{t_2}^{(zu)_2}(z_2) \cdot ... \cdot q_{t_{n-1}}^{(zu)_{n-1}}(z_{n-1}) q_{t_n}^{(zu)_n}(\alpha)}{q_{t_1}^{(zu)_1}(u_1) q_{t_2}^{(zu)_2}(u_2) \cdot ... \cdot q_{t_{n-1}}^{(zu)_{n-1}}(u_{n-1}) q_{t_n}^{(zu)_n}(u_n)} \right)^{-1} \\
& = & \dfrac{q_{t_n}^{(zu)_n}(x)}{\sum \limits_{\alpha \in X^{\{t\}}} q_{t_n}^{(zu)_n}(\alpha)} = q_t^z(x).
\end{eqnarray*}
\end{proof}

Note that the formula~\eqref{Q1->P} can be written in the following equivalent forms
\begin{equation}
\label{Q1->P-2}
P_\Lambda(x) = \left( \sum \limits_{\alpha \in X^\Lambda} \prod \limits_{j=1}^n \frac{q_{t_j}^{(\alpha x)_j} (\alpha_j)}{q_{t_j}^{(\alpha x)_j} (x_j)} \right)^{-1} = \frac{q_{t_1}^{x_{\Lambda \backslash \{t_1\}}} (x_1)}{ \sum \limits_{\alpha \in X^{\Lambda \backslash \{t_n\}}} \frac{q_{t_1}^{x_{\Lambda \backslash \{t_1\}}} (\alpha_1)}{q_{t_n}^{\alpha} (x_n)}  \prod \limits_{j=2}^{n-1} \frac{q_{t_j}^{(\alpha x)_j} (\alpha_j)}{q_{t_j}^{(\alpha x)_j} (x_j)}}
\end{equation}
for any $x \in X^\Lambda$ and any enumeration $\Lambda = \{t_1, t_2, ..., t_n\}$ of points in $\Lambda$, $n = \vert \Lambda \vert$. Such representations make it obvious that the right-hand side of~\eqref{Q1->P} does not depend on the choice of configuration $u \in X^\Lambda$.

\subsection{Weakly positive random fields}

Now let us consider random fields $P_\Lambda$ which can take zero values.

A random field $P_\Lambda$ we call \emph{weakly positive} if for each $t \in \Lambda$ there exists an element $\theta_t \in X^t$ such that $P_\Lambda(\theta_t z)>0$ for any $z \in X^{\Lambda \backslash \{t\}}$.

Thus, a weakly positive random field $P_\Lambda$ takes strictly positive values on those configurations $x \in X_\Lambda$ which have at least one component $x_t$ coinciding with $\theta_t$, $t \in \Lambda$. It is clear, that for such random field, all the one-point conditional probabilities~\eqref{P->Q1} are defined, and  $Q_t^z(\theta_t) > 0$ for any boundary condition $z \in X^{\Lambda \backslash \{t\}}$, $t \in \Lambda$.

%In this case,
%$$
%\sum \limits_{u \in X^t} P_\Lambda(uz) \ge P_\Lambda(\theta_t z) > 0
%$$
%for all $t \in \Lambda$ and $z \in X^{\Lambda \backslash \{t\}}$, and hence, for the weakly positive random field, all the one-point conditional probabilities~\eqref{P->Q1} are defined. It is clear that $Q_t^z(\theta_t) > 0$ for any boundary condition $z \in X^{\Lambda \backslash \{t\}}$, $t \in \Lambda$.

A family $\boldsymbol{\rm Q}_1 = \{q_t^z, z \in X^{\Lambda \backslash \{t\}}, t \in \Lambda \}$ of one-point probability distributions $q_t^z$ on $X^t$ parameterized by boundary conditions $z \in X^{\Lambda \backslash \{t\}}$, $t \in \Lambda$, will be called \emph{weakly positive} if for each $t \in \Lambda$ there exists an element $\theta_t \in X^t$ for which $q_t^z(\theta_t)>0$ under any boundary condition $z \in X^{\Lambda \backslash \{t\}}$. Such an element $\theta_t$ we will call \emph{positivity point} of the system $\boldsymbol{\rm Q}_1$ at $t \in \Lambda$ (compare with the analogous notion in~\cite{DN04}).

Note that if for a weakly positive family  $\boldsymbol{\rm Q}_1 = \{q_t^z, z \in X^{\Lambda \backslash \{t\}}, t \in \Lambda \}$ of probability distributions, all its positivity points $\theta_t$, $t \in \Lambda$, coincide with some specified element $\theta \in X$, then $\theta$ is called the vacuum and $\boldsymbol{\rm Q}_1$ is called \emph{vacuum system}. Such systems play important role in the mathematical statistical physics.

The following statement holds true.

\begin{T}
\label{Th-Q1->P}
For a given weakly positive system $\boldsymbol{\rm Q}_1 = \{q_t^z, z \in X^{\Lambda \backslash \{t\}}, t \in \Lambda \}$, there exists a unique weakly positive random field $P_\Lambda$ compatible with it if and only if the elements of $\boldsymbol{\rm Q}_1$ satisfy the following consistency conditions: for all $t,s \in \Lambda$ and $x \in X^t$, $y \in X^s$, $z \in X^{\Lambda \backslash \{t,s\}}$, it holds
$$
q_t^{zy} (x) q_s^{zx} (\theta_s) q_t^{z\theta_s} (\theta_t) q_s^{z\theta_t} (y) = q_s^{zx} (y) q_t^{zy} (\theta_t) q_s^{z\theta_t} (\theta_s) q_t^{z\theta_s} (x),
$$
where $\theta_t$ and $\theta_s$ are positivity points of the system $\boldsymbol{\rm Q}_1$.
\end{T}
\begin{proof}
The necessity is obvious. Let us prove the sufficiency.

For any $x \in X^\Lambda$, put
$$
P_\Lambda(x) = Z^{-1}(\theta_{t_1} \theta_{t_2} ... \theta_{t_{\vert \Lambda \vert}}) \frac{q_{t_1}^{\theta_{t_2} \theta_{t_3} ... \theta_{t_{\vert \Lambda \vert}}}(x_{t_1}) q_{t_2}^{x_{t_1} \theta_{t_3} ... \theta_{t_{\vert \Lambda \vert}}}(x_{t_2}) \cdot ... \cdot q_{t_{\vert \Lambda \vert}}^{x_{t_1} x_{t_2} ... x_{t_{\vert \Lambda \vert -1}}}(x_{t_{\vert \Lambda \vert}})}{q_{t_1}^{\theta_{t_2} \theta_{t_3} ... \theta_{t_{\vert \Lambda \vert}}}(\theta_{t_1}) q_{t_2}^{x_{t_1} \theta_{t_3} ... \theta_{t_{\vert \Lambda \vert}}}(\theta_{t_2}) \cdot ... \cdot q_{t_{\vert \Lambda \vert}}^{x_{t_1} x_{t_2} ... x_{t_{\vert \Lambda \vert -1}}}(\theta_{t_{\vert \Lambda \vert}})},
$$
where
$$
Z(\theta_{t_1} \theta_{t_2} ... \theta_{t_{\vert \Lambda \vert}}) = \sum \limits_{\alpha \in X^\Lambda} \frac{q_{t_1}^{\theta_{t_2} \theta_{t_3} ... \theta_{t_{\vert \Lambda \vert}}}(\alpha_{t_1}) q_{t_2}^{\alpha_{t_1} \theta_{t_3} ... \theta_{t_{\vert \Lambda \vert}}}(\alpha_{t_2}) \cdot ... \cdot q_{t_{\vert \Lambda \vert}}^{\alpha_{t_1} \alpha_{t_2} ... \alpha_{t_{\vert \Lambda \vert -1}}}(\alpha_{t_{\vert \Lambda \vert}})}{q_{t_1}^{\theta_{t_2} \theta_{t_3} ... \theta_{t_{\vert \Lambda \vert}}}(\theta_{t_1}) q_{t_2}^{\alpha_{t_1} \theta_{t_3} ... \theta_{t_{\vert \Lambda \vert}}}(\theta_{t_2}) \cdot ... \cdot q_{t_{\vert \Lambda \vert}}^{\alpha_{t_1} \alpha_{t_2} ... \alpha_{t_{\vert \Lambda \vert -1}}}(\theta_{t_{\vert \Lambda \vert}})},
$$
$\theta_{t_j}$ is a positivity point of $\boldsymbol{\rm Q}_1$ at $t_j$, $1 \le j \le \vert \Lambda \vert$, and $\Lambda = \{t_1, t_2, ..., t_{\vert \Lambda \vert}\}$ is some enumeration of the elements of  $\Lambda$.

Using the same reasoning as in the proof of Theorem~\ref{Th-Q1->P+}, one can show that the probability distribution $P_\Lambda$ does not depend on the enumeration of the elements in $\Lambda$ and on the choice of a positivity point at each $t \in \Lambda$ (if there are several positivity points at $t$). Finally, it is not difficult to verify that $\boldsymbol{Q}_1(P_\Lambda) = \boldsymbol{\rm Q}_1$.
\end{proof}

\subsection{Notes on previous results}

A large number of works are devoted to the question of describing finite random fields by consistent systems of conditional distributions. Let us note only those which results are close to obtained in this paper.

Apparently, the first necessary and sufficient conditions for the existence and uniqueness of a two-dimensional ($\vert \Lambda \vert=2$) joint distribution for a given system of one-point probability distributions $\boldsymbol{\rm Q}_1$ were obtained by Abrahams and Thomas~\cite{AbTh}. It should be noted that their conditions are not expressed in terms of the elements of the given system $\boldsymbol{\rm Q}_1$. At the same time, from the point of view of the statement of the problem, it is natural to impose conditions directly on the elements of $\boldsymbol{\rm Q}_1$.

Abrahams and Thomas' conditions were considered in several subsequent works (see, for example, Arnold and Press~\cite{ArPr} and Berti et al.~\cite{Berti}). Arnold and Press~\cite{ArPr} introduced the equivalent conditions expressed in terms of the elements of the given system $\boldsymbol{\rm Q}_1$ in the cases $\vert \Lambda \vert=2$ and $\vert \Lambda \vert=3$. These conditions coincides with our consistency conditions in the considered by them cases. Nevertheless, none of the mentioned authors provide an explicit formula expressing $P_\Lambda$ in terms of the elements of $\boldsymbol{\rm Q}_1$.

%Concerning our results, for any finite set $\Lambda$, we obtained (Theorem~\ref{Th-Q1->P+}) necessary and sufficient conditions for the existence and uniqueness of a finite random field $P_\Lambda$ with a given system of one-point conditional distributions $\boldsymbol{\rm Q}_1$. These conditions have a rather simple form and are expressed directly in terms of the elements of $\boldsymbol{\rm Q}_1$. Moreover, we suggest an (explicit) formula (see formula~\eqref{Q1->P}) connecting $P_\Lambda$ with the elements of $\boldsymbol{\rm Q}_1$ for all finite sets $\Lambda$.

Such a formula was obtained by Chen~\cite{Chen}. This formula, after the necessary simplifications, coincides with our formula~\eqref{Q1->P}. However, the consistency conditions used by Chen are, in fact, a requirement for the correctness of this formula, while the validation criterion of correctness is not provided.
%The consistency conditions~\eqref{Q1-sogl} we used is such a criterion.
Note also that Chen does not discuss the uniqueness of the constructed random field $P_\Lambda$.

Gurevich~\cite{Gurev} gives necessary and sufficient conditions for the existence of a finite random field with a given system $\boldsymbol{\rm Q}_1$. He proposed the formula connecting $P_\Lambda$ with the elements of $\boldsymbol{\rm Q}_1$. However, like in the work by Chen~\cite{Chen}, the conditions proposed by Gurevich, in their essence, are the requirement of the correctness of this formula.

Regarding the uniqueness of a random field with a given system of conditional probabilities, we note that the answer was given mainly for positive systems of consistent distributions. For systems admitting zeros, the uniqueness conditions are given by Gurevich~\cite{Gurev} (see also Arnold and Press~\cite{ArPr}).

%It is not difficult to see that Theorem~\ref{Th-Q1->P} remains true for vacuum systems as well.

\section{Gibbs representation of random fields}
\label{Sec-Gibbs}

In many practical situations, the statistical data are often presented in the form of a system of conditional probability distributions. In particular, this is the case in mathematical statistical physics, where its main object, Gibbs random field, is described by means of a special system of conditional distributions, the Gibbs specification. The Gibbsian form of the conditional distribution of a random field is universal, namely, all systems of consistent conditional distributions can be represented in the Gibbs form (see~\cite{DN19}).

The Gibbs representation makes it very easy to construct specific examples of systems of consistent probability distributions. Moreover, the Gibbs representation allows one to describe various classes of random fields, such as, for example, the classes of Markov, stationary, and weakly dependent random fields.

The problem of Gibbsianess of a random field $P_\Lambda$ was considered by Griffeath~\cite{KSK}. He showed that any positive distribution in a finite volume can be represented in the Gibbs form~\eqref{P-Gibbs} with a vacuum potential without long-range interactions. In this case, an explicit form of the potential is written out in terms of $P_\Lambda$ using the M\"{o}bius formula.

Note that the Gibbs representation is especially suitable when the interaction potential has a simple form, while the application of the M\"{o}bius formula leads to a very complex expression for it and may have only theoretical interest. In addition, as noted by Gandolfi and Lenarda~\cite{GL}, when applying the M\"{o}bius formula, the explicit dependence of the potential on the values of the spins of physical systems is lost.

In this section, using analogues of the methods proposed in~\cite{DN19}, we present (based on the energy conservation law) the definition of a Gibbs finite random field in a much more general form (without using of the notion of potential) than the commonly accepted one. We also show that any finite positive random field is Gibbsian. Applied to Markov random fields, the proposed approach leads to a generalized version of the Hammersley--Clifford theorem.

In what follows, we will consider positive finite random fields only.

\subsection{Gibbs random fields with a given potential}

Note (see also~\cite{DN19}) that a function $P_\Lambda(x)$, $x \in X^\Lambda$, is a probability distribution on $X^\Lambda$ if and only if it is representable in the form
\begin{equation}
\label{P-Gibbs}
P_\Lambda(x) = \frac{\exp\{-H_\Lambda(x)\}}{\sum \limits_{\alpha \in X^\Lambda} \exp\{-H_\Lambda(\alpha)\}}, \qquad x \in X^\Lambda,
\end{equation}
where $H_\Lambda$ is some function on $X^\Lambda$. It is clear that this formula is unlikely to be useful if one does not impose certain restrictions on $H_\Lambda$.

In statistical physics, it is assumed that the function $H_\Lambda$ should be a Hamiltonian (potential energy), and it is postulated that the Hamiltonian should have the following form
\begin{equation}
\label{H-Phi}
H_\Lambda(x) = \sum \limits_{\emptyset \neq V \subseteq \Lambda} \Phi_V(x_V), \qquad x \in X^\Lambda,
\end{equation}
where interaction potential $\Phi$ is a family $\{\Phi_V, V \subseteq \Lambda\}$ of functions $\Phi_V$ on $X^V$, $V \subseteq \Lambda$. This approach makes it possible to build models of  physical systems based on potential.

In the theory of finite random fields, by a \emph{Gibbs random field with a given potential $\Phi$} it is accepted to call (see, for example,~\cite{KSK}) a random field $P_\Lambda$, which is representable in the form~\eqref{P-Gibbs} with the Hamiltonian $H_\Lambda$ constructed by the potential $\Phi$.

In the framework of the problem of description of random fields by systems of conditional distributions, we give an equivalent definition of a Gibbs random field $P_\Lambda$ (with a potential $\Phi$) in terms of $\boldsymbol{Q}_1(P_\Lambda)$, which is in accordance with the one used in the case of infinite volumes.

For a given potential $\Phi = \{\Phi_V, V \subseteq \Lambda\}$, consider the Gibbs system $\boldsymbol{\rm Q}_1^\Phi = \{q_t^z, z \in X^{\Lambda \backslash \{t\}}, t \in \Lambda \}$ of probability distributions
\begin{equation}
\label{QPhi}
q_t^z(x) = \frac{\exp\{-H_t^z(x)\}}{\sum \limits_{\alpha \in X^t} \exp\{-H_t^z(\alpha)\}}, \qquad x \in X^t,
\end{equation}
where
\begin{equation}
\label{H-tz-Phi}
H_t^z(x) = \sum \limits_{J \subseteq \Lambda \backslash \{t\}} \Phi_{t \cup J}(x z_J), \qquad x \in X^t, z \in X^{\Lambda \backslash \{t\}}, t \in \Lambda.
\end{equation}
Note that if $H_\Lambda$ is defined by~\eqref{H-Phi}, then
$$
H_t^z(x) = H_\Lambda(xz) - H_{\Lambda \backslash \{t\}}(z), \qquad x \in X^t, z \in X^{\Lambda \backslash \{t\}}, t \in \Lambda.
$$
It is not difficult to see, that the elements of $\boldsymbol{\rm Q}_1^\Phi$ satisfy the consistency conditions~\eqref{Q1-sogl}. Therefore, according to Theorem~\ref{Th-Q1->P+}, $\boldsymbol{\rm Q}_1^\Phi$ uniquely determines a random field $P_\Lambda$ such that $\boldsymbol{Q}_1(P_\Lambda) = \boldsymbol{\rm Q}_1^\Phi$. Hence, we come to the following definition of the Gibbs random field with a potential $\Phi$ formulated in terms of the elements of $\boldsymbol{Q}_1(P_\Lambda)$ only.

A random field $P_\Lambda$ will be called a \emph{Gibbs random field with potential $\Phi = \{\Phi_V, V \subseteq \Lambda\}$} if $\boldsymbol{Q}_1(P_\Lambda) = \boldsymbol{\rm Q}_1^\Phi$ where $\boldsymbol{\rm Q}_1^\Phi$ is the Gibbs system constructed by the potential $\Phi$.

%the elements of its one-point conditional distribution $\boldsymbol{Q}_1(P_\Lambda) = \{Q_t^z, z \in X^{\Lambda \backslash \{t\}}, t \in \Lambda \}$ are representable in the form
%\begin{equation}
%\label{Q(P)-H-Phi}
%Q_t^z(x) = \frac{\exp\{-H_t^z(x)\}}{\sum \limits_{\alpha \in X^\Lambda} \exp\{-H_t^z(\alpha)\}}, \qquad x \in X^t.
%\end{equation}

It is clear that if $P_\Lambda$ is a Gibbs random field with potential $\Phi$ in the conventional sense, then its one point conditional distributions have the Gibbs form~\eqref{QPhi}. Conversely, let the elements of the system $\boldsymbol{Q}_1(P_\Lambda)$ have the form~\eqref{QPhi}. According to the formula~\eqref{Q1->P-2}, we can write
\begin{eqnarray*}
P_\Lambda(x) & = & \left( \sum \limits_{\alpha \in X^\Lambda} \prod \limits_{j=1}^n \dfrac{Q_{t_j}^{(\alpha x)_j} (\alpha_j)}{Q_{t_j}^{(\alpha x)_j} (x_j)} \right)^{-1} = \left( \sum \limits_{\alpha \in X^\Lambda} \prod \limits_{j=1}^n \dfrac{\exp\{ - H_{t_j}^{(\alpha x)_j} (\alpha_j)\}}{\exp\{ - H_{t_j}^{(\alpha x)_j} (x_j)\}} \right)^{-1} \\
\\
& = & \left( \sum \limits_{\alpha \in X^\Lambda} \exp \left\{ \sum \limits_{j=1}^n \left( H_{t_j}^{(\alpha x)_j} (x_j) - H_{t_j}^{(\alpha x)_j} (\alpha_j)\right) \right\} \right)^{-1},
\end{eqnarray*}
where $\Lambda = \{t_1, t_2, ..., t_n\}$ is some enumeration of the points in $\Lambda$, $n = \vert \Lambda \vert$. Since
$$
%\begin{array}{l}
\sum \limits_{j=1}^n \left( H_{t_j}^{(\alpha x)_j}(x_j) - H_{t_j}^{(\alpha x)_j}(\alpha_j) \right) %\\
%\\
% \qquad = \sum \limits_{j=1}^n \left( H_\Lambda(\alpha_1 ... \alpha_{j-1} x_j x_{j+1} ... x_n)  - H_\Lambda(\alpha_1 ... \alpha_{j-1} \alpha_j x_{j+1} ... x_n) \right)
= H_\Lambda(x) - H_\Lambda(\alpha),
%\end{array}
$$
finally we obtain
$$
P_\Lambda(x) = \left( \sum \limits_{\alpha \in X^\Lambda} \exp \{H_\Lambda(x) - H_\Lambda(\alpha)\} \right)^{-1} = \frac{\exp\{-H_\Lambda(x)\}}{\sum \limits_{\alpha \in X^\Lambda} \exp\{-H_\Lambda(\alpha)\}}, \qquad x \in X^\Lambda.
$$

Given the above definitions of a Gibbs random field are essentially based on the fact that the Hamiltonian is the sum of the interaction potentials. Such a representation of the Hamiltonian cannot be considered as its definition since in the case of this approach, the properties of the Hamiltonian are not postulated but derived from the properties of the potential. In~\cite{DN19}, based on the energy conservation law, an axiomatic definition of the Hamiltonian was given by means of its intrinsic properties. This definition is based on the concept of a transition energy field introduced in~\cite{DN19}.

Applying the approach developed in~\cite{DN19} to finite random fields, below we give the definition of Gibbs random field without using the notion of interaction potential. Note that some of our arguments literally repeat those used in the case of infinite volumes. This is done for the sake of completeness.

\subsection{Axiomatic definition of Hamiltonian and Gibbs random fields}

First of all, we give the following (axiomatic) definition of the Hamiltonian (in finite volume) which does not use the notion of potential (see~\cite{DN19}).

A family $H_1 = \{ H_t^z, z \in X^{\Lambda \backslash \{t\}}, t \in \Lambda\}$ of functions $H_t^z$ on $X^t$ parameterized by boundary conditions $z \in X^{\Lambda \backslash \{t\}}$, $t \in \Lambda$, we call (\emph{one-point}) \emph{Hamiltonian} in $\Lambda$ if its elements satisfy the following conditions: for all $t,s \in \Lambda$ and $x,u \in X^t$, $y,v \in X^s$, $z \in X^{\Lambda \backslash \{t,s\}}$, it holds
\begin{equation}
\label{H-sogl}
H_t^{zy}(x) + H_s^{zx}(v) + H_t^{zv}(u) + H_s^{zu}(y) = H_t^{zy}(u) + H_s^{zu}(v) + H_t^{zv}(x) + H_s^{zx}(y).
\end{equation}

The basic example of a one-point Hamiltonian is a Hamiltonian $H_1^\Phi = \{ H_t^z, z \in X^{\Lambda \backslash \{t\}}, t \in \Lambda\}$ constructed by the potential $\Phi = \{\Phi_V, V \subseteq \Lambda\}$ (see formula~\eqref{H-tz-Phi}). Another example is the Hamiltonian $H_1^\theta (P_\Lambda) = \{ H_t^z, z \in X^{\Lambda \backslash \{t\}}, t \in \Lambda\}$ constructed by a random field $P_\Lambda$ and a fixed configuration $\theta \in X^\Lambda$ according to the formula
$$
H_t^z(x) = \ln \frac{P_\Lambda(xz)}{P_\Lambda(\theta_t z)}, \qquad x \in X^t, z \in X^{\Lambda \backslash \{t\}}, t \in \Lambda.
$$
It is easy to check that the elements both of $H_1^\Phi$ and $H_1^\theta(P_\Lambda)$ satisfy the conditions~\eqref{H-sogl}.

Consistency conditions~\eqref{H-sogl} have a very simple physical interpretation. First note that they can be written in the following form
$$
H_t^{zy}(u) - H_t^{zy}(x) + H_s^{zu}(v) - H_s^{zu}(y) = H_s^{zx}(v) - H_s^{zx}(y) + H_t^{zv}(u) - H_t^{zv}(x)
$$
For each $t \in \Lambda$, consider the function
\begin{equation}
\label{Delta-H}
\Delta_t^z(x,u) = H_t^z(u) - H_t^z(x), \qquad x,u \in X^t, z \in X^{\Lambda \backslash \{t\}},
\end{equation}
which represents an amount of energy needed under fixed boundary condition $z$ to change the state of the system at point $t$ from $x$ to $u$. It is clear that for each $t \in \Lambda$ and $z \in X^{\Lambda \backslash \{t\}}$, the function $\Delta_t^z$ satisfy the relation
\begin{equation}
\label{Delta1Lambda-1}
\Delta_t^z(x,u) = \Delta_t^z(x,\alpha) + \Delta_t^z(\alpha,u), \qquad x,u,\alpha \in X^t,
\end{equation}
which is in full compliance with the energy conservation law. Conditions~\eqref{H-sogl} in terms of $\Delta_t^z$ take the following form: for all $t,s \in \Lambda$ and $x,u \in X^t$, $y,v \in X^s$, $z \in X^{\Lambda \backslash \{t,s\}}$, it holds
\begin{equation}
\label{Delta1Lambda-2}
\Delta_t^{zy}(x,u) + \Delta_s^{zu}(y,v) = \Delta_s^{zx}(y,v) + \Delta_t^{zv}(x,u).
\end{equation}
This, in turn, can be interpreted as follows. Suppose it is necessary to move a physical system from the state $xyz$ to the state $uvz$, where $x,u \in X^t$, $y,v \in X^s$, $z \in X^{\Lambda \backslash \{t,s\}}$, $t,s \in \Lambda$. This can be done in two ways. First, change the state of the system at the point $t$ from $x$ to $u$ under the boundary condition $zy$, and then at the point $s$ from $y$ to $v$ already under the boundary condition $zu$. The energy expended in this case is $\Delta_t^{zy}(x,u) + \Delta_s^{zu}(y,v)$. One can proceed also as follows: first, starting from the point $s$, change the state from $y$ to $v$ under the boundary condition $zx$, and then, under the boundary condition $zv$, change it at the point $t$ from $x$ to $u$ expending the energy $\Delta_s^{zx}(y,v) + \Delta_t^{zv}(x,u)$. Naturally, the same amount of energy must be spent in both cases.

Note also, that from the physical point of view, it is more natural to consider the difference of energies instead of the energies themselves. Thus, following~\cite{DN19}, we introduce the notion of transition energy field in finite volume.

A family $\boldsymbol{\Delta}_1 = \{ \Delta_t^z, z \in X^{\Lambda \backslash \{t\}}, t \in \Lambda\}$ of functions $\Delta_t^z$ on pairs of configurations from $X^t$ and parameterized by boundary conditions $z \in X^{\Lambda \backslash \{t\}}$, $t \in \Lambda$, we call a (\emph{one-point}) \emph{transition energy field} in $\Lambda$ if its elements satisfy the consistency conditions~\eqref{Delta1Lambda-1} and~\eqref{Delta1Lambda-2}.

The fundamental connection between the consistent system $\boldsymbol{\rm Q}_1$ and the transition energy field $\boldsymbol{\Delta}_1$ is established in the theorem below. Since the proof is similar to the proof of Theorem 2.3 in~\cite{DN19}, we give it briefly.

\begin{T}
\label{Th-Q1-Delta1}
A family of functions $\boldsymbol{\rm Q}_1 = \{q_t^z, z \in X^{\Lambda \backslash \{t\}}, t \in \Lambda\}$ forms a system of positive one-point probability distributions consistent in the sense~\eqref{Q1-sogl} if and only if its elements have a Gibbs form
\begin{equation}
\label{Q-Gibbs-Delta}
q_t^z(x) = \frac{\exp\{\Delta_t^z(x,u)\}}{\sum \limits_{\alpha \in X^t} \exp\{\Delta_t^z(\alpha,u)\}}, \qquad x \in X^t,
\end{equation}
where $u \in X^t$ and $\boldsymbol{\Delta}_1 = \{\Delta_t^z, z \in X^{\Lambda \backslash \{t\}}, t \in \Lambda\}$ is a one-point transition energy field. Herewith, the transition energy field $\boldsymbol{\Delta}_1^\Lambda$ is uniquely determined by $\boldsymbol{\rm Q}_1$.
\end{T}
\begin{proof}
Let the elements of the system $\boldsymbol{\rm Q}_1 = \{q_t^z, z \in X^{\Lambda \backslash \{t\}}, t \in \Lambda\}$ of positive probability distributions satisfy the consistency conditions~\eqref{Q1-sogl}. For each $t \in \Lambda$ and $z \in X^{\Lambda \backslash \{t\}}$, put
\begin{equation}
\label{DeltaQ}
\Delta_t^z(x,u) = \ln \frac{q_t^z(x)}{q_t^z(u)}, \qquad x,u \in X^t.
\end{equation}
It is easy to verify that the elements of the system $\boldsymbol{\Delta}_1 = \{\Delta_t^z, z \in X^{\Lambda \backslash \{t\}}, t \in \Lambda\}$ satisfy the conditions~\eqref{Delta1Lambda-1} and~\eqref{Delta1Lambda-2} , and thus, $\boldsymbol{\Delta}_1^\Lambda$ is a transition energy field.

%Indeed, for any $x,u, \alpha \in X^t$, we have
%$$
%\Delta_t^z(x,u) = \ln \frac{q_t^z(x) q_t^z(\alpha)}{q_t^z(u) q_t^z(\alpha)} = \Delta_t^z (x,\alpha) + \Delta_t^z (\alpha,u).
%$$
%Further, by the consistency conditions of the elements of $\boldsymbol{\rm Q}_1$, for all $t,s \in \Lambda$ and $x,u \in X^t$, $y,v \in X^s$, $z \in X^{\Lambda \backslash \{ t,s\}}$, we can write
%$$
%\Delta_t^{zy} (x,u) + \Delta_s^{zu} (y,v) =  \ln \dfrac{q_t^{zy} (x)}{q_t^{zy} (u)} \cdot \dfrac{q_s^{zu} (y)}{q_s^{zu} (v)} = \ln \dfrac{q_s^{zx} (y)}{q_s^{zx} (v)} \cdot \dfrac{q_t^{zv} (x)}{q_t^{zv} (u)} = \Delta_s^{zx} (y,v) + \Delta_t^{zv} (x,u).
%$$

Let now $\boldsymbol{\Delta}_1 = \{\Delta_t^z, z \in X^{\Lambda \backslash \{t\}}, t \in \Lambda\}$ be a transition energy field. For each $t \in \Lambda$ and $z \in X^{\Lambda \backslash \{t\}}$, put
$$
q_t^z(x) = \frac{\exp \{\Delta_t^z(x,u)\}}{\sum \limits_{\alpha \in X} {\exp \{\Delta_t^z (\alpha,u)\} }}, \qquad x,u \in X^t.
$$
First, note that due to~\eqref{Delta1Lambda-1}, the values of $q_t^z$ do not depend on the choice of $u$.
%since by for any $w \in X^t$, we have
%$$
%\frac{\exp \{\Delta_t^z(x,u)\}}{\sum \limits_{\alpha \in X} {\exp \{\Delta_t^z(\alpha,u)\}}} = \frac{\exp \{\Delta_t^z(x,w) + \Delta_t^z(w,u)\}}{\sum \limits_{\alpha \in X} {\exp \{\Delta_t^z(\alpha,w) + \Delta_t^z(w,u)\} }} = \frac{\exp \{\Delta_t^z(x,w)\}}{\sum \limits_{\alpha \in X} {\exp \{\Delta_t^z(\alpha,w)\}}}.
%$$
Let us verify that the elements of the family $\left\{ q_t^z, z \in X^{\Lambda \backslash \{t\}}, t \in \Lambda \right\}$ satisfy the consistency conditions~\eqref{Q1-sogl}. This statement follows from the chain of equalities
$$
\ln \dfrac{q_t^{zy} (x)}{q_t^{zy} (u)} \cdot \dfrac{q_s^{zu} (y)}{q_s^{zu} (v)} = \Delta_t^{z y} (x,u) + \Delta_s^{zu} (y,v) = \Delta_s^{zx} (y,v) + \Delta_t^{zv} (x,u) = \ln \dfrac{q_s^{zx} (y)}{q_s^{zx} (v)} \cdot \dfrac{q_t^{zv} (x)}{q_t^{zv} (u)},
$$
$x,u \in X^t$, $y,v \in X^s$,  $z \in X^{\Lambda \backslash \{ t,s\}}$, $t,s \in \Lambda$, which is true by virtue of the consistency conditions~\eqref{Delta1Lambda-1} of the elements of $\boldsymbol{\Delta}_1^\Lambda$.
%
%To do so, it is sufficient to show that for all  it holds
%$$
%\frac{q_t^{zy} (x)}{q_t^{zy} (u)} \cdot \frac{q_s^{zu} (y)}{q_s^{zu} (v)} = \frac{q_s^{zx} (y)}{q_s^{zx} (v)} \cdot \frac{q_t^{zv} (x)}{q_t^{zv} (u)}.
%$$
%The validity of the last relation
\end{proof}

We propose the following definition of a Gibbs finite random field, based on the introduced notion of the transition energy field.

A random field $P_\Lambda$ we call \emph{Gibbsian} if the elements of $\boldsymbol{Q}_1(P_\Lambda)$ can be represented in the Gibbs form
$$
Q_t^z(x) = \frac{\exp\{\Delta_t^z(x,u)\}}{\sum \limits_{\alpha \in X^t} \exp\{\Delta_t^z(\alpha,u)\}}, \qquad x \in X^t,
$$
where $u \in X^t$ and $\boldsymbol{\Delta}_1 = \{ \Delta_t^z, z \in X^{\Lambda \backslash \{t\}}, t \in \Lambda\}$ is the transition energy field.

It is easy to see that the definition of a Gibbs random field based on the notion of transition energy field includes all the definitions given above.

The direct corollary of Theorem~\ref{Th-Q1-Delta1} is the following result.

\begin{T}
\label{Th-allP-Gibbs}
Any finite random field is Gibbsian.
\end{T}
\begin{proof}
Since the elements of the one-point conditional distribution $\boldsymbol{Q}_1(P_\Lambda)$ of any random field $P_\Lambda$ satisfy the consistency conditions~\eqref{Q1-sogl}, by Theorem~\ref{Th-Q1-Delta1} they have the Gibbs form~\eqref{Q-Gibbs-Delta} with uniquely determined transition energy field $\boldsymbol{\Delta}_1 = \{\Delta_t^z, z \in X^{\Lambda \backslash \{t\}}, t \in \Lambda\}$. Herewith,
$$
\Delta_t^z(x,u) = \ln \frac{Q_t^z(x)}{Q_t^z(u)} = \ln \frac{P_\Lambda(xz)}{P_\Lambda(uz)}, \qquad x,u \in X^t.
$$
\end{proof}

Thus, there is a one-to-one correspondence between the positive random field $P_\Lambda$ and the one-point transition energy field $\boldsymbol{\Delta}_1$.

From Theorem~\ref{Q-Gibbs-Delta}, it directly follows the validity of the next result.

\begin{T}
\label{Th-QDLambda-Gibbs}
The family of functions $\boldsymbol{\rm Q}_1 = \{q_t^z, z \in X^{\Lambda \backslash \{t\}}, t \in \Lambda\}$ forms a system of positive one-point probability distributions satisfying the consistency conditions~\eqref{Q1-sogl} if and only if its elements have a Gibbs form
$$
q_t^z(x) = \frac{\exp\{ -H_t^z(x)\}}{\sum \limits_{\alpha \in X^t} \exp\{ -H_t^z(\alpha)\}}, \qquad x \in X^t, z \in X^{\Lambda \backslash \{t\}}, t \in \Lambda,
$$
where $H_1 = \{ H_t^z, z \in X^{\Lambda \backslash \{t\}}, t \in \Lambda\}$ is some one-point Hamiltonian.
\end{T}

Unlike the transition energy field, the Hamiltonian $H_1^\Lambda$ for a consistent system $\boldsymbol{\rm Q}_1$ (and, therefore, compatible with it random field) is determined not uniquely but up to an additive function on boundary conditions. Nevertheless, all such Hamiltonians lead to the same consistent system $\boldsymbol{\rm Q}_1$. This remark allows us to give the following definition of a Gibbs random field in a finite volume.

A random field $P_\Lambda$ is called \emph{Gibbsian} if the elements of its one-point conditional distribution $\boldsymbol{Q}_1(P_\Lambda)$ are representable in a Gibbs form
$$
Q_t^z(x) = \frac{\exp\{-H_t^z(x)\}}{\sum \limits_{\alpha \in X^t} \exp\{-H_t^z(\alpha)\}}, \qquad x \in X^t,
$$
where $H_1 = \{ H_t^z, z \in X^{\Lambda \backslash \{t\}}, t \in \Lambda\}$ is some one-point Hamiltonian.

\subsection{Markov finite random fields}

Most often, the question of Gibbsianess was considered for Markov finite random fields. The main result here is considered to be the Hammersley--Clifford theorem on the equivalence of Markov and Gibbs random fields with a pair interaction potential of the nearest neighbors, first presented in an unpublished paper~\cite{HammCliff}. The proof of this theorem has been repeatedly presented under various conditions both on $\Lambda$ and on the state space $X$. Let us mention, for example, works by Grimmet~\cite{Grimmet}, Preston~\cite{Preston}, Spitzer~\cite {Spitzer}, and Besag~\cite{Besag}. However, as it will be shown below, the Hammersley--Clifford theorem is a direct consequence of a much more general theorem (see Theorem~\ref{Th-Markov-Gibbs}).

The definition of Markov random fields is based on the concept of a neighborhood system. A family $\partial = \{\partial t, t \in \Lambda\}$ of sets $\partial t \subset \Lambda$ is called a neighborhood system in $\Lambda$ if $t \notin \partial t$ and $s \in \partial t$ if and only if $t \in \partial s$, $t,s \in \Lambda$.

A random field $P_\Lambda$ is called \emph{Markov random field} (with respect to $\partial$) if its conditional distributions satisfy the following property: for any $t \in \Lambda$ and $z \in X^{\Lambda \backslash \{t\}}$, it holds
$$
Q_t^z(x) = Q_t^{z_{\partial t}}(x), \qquad x \in X^t.
$$

It is not difficult to verify that a random field $P_\Lambda$ will be a Markov random field if and only if for each $t \in \Lambda$ and any $x,u \in X^t$, $z \in X^{\partial t}$, $y,v \in X^{\Lambda \backslash (t \cup \partial t)}$, it holds
$$
Q_t^{zy}(x) = Q_t^{zv}(x).
$$

Let us now present conditions on the transition energy field, under which the corresponding finite random field is Markov.

\begin{prop}
\label{Prop-Markov-Delta}
A random field $P_\Lambda$ is a Markov random field if and only if the elements of corresponding transition energy field $\boldsymbol{\Delta}_1 = \{ \Delta_t^z, z \in X^{\Lambda \backslash \{t\}}, t \in \Lambda\}$ satisfy the following conditions: for all $t \in \Lambda$ and any $x,u \in X^t$, $z \in X^{\partial t}$, $y,v \in X^{\Lambda \backslash (t \cup \partial t)}$, it holds
\begin{equation}
\label{Delta-Markov}
\Delta_t^{zy}(x,u) = \Delta_t^{zv}(x,u), \qquad x,u \in X^t.
\end{equation}
\end{prop}
\begin{proof}
Let $\boldsymbol{Q}_1(P_\Lambda)$ be a system of one-point conditional distributions of the Markov random field $P_\Lambda$. Then for all $t \in \Lambda$ and $x,u \in X^t$, $z \in X^{\partial t}$, $y,v \in X^{\Lambda \backslash (t \cup \partial t)}$, we have
$$
\Delta_t^{zy}(x,u) = \ln \frac{Q_t^{zy}(x)}{Q_t^{zy}(u)} = \ln \frac{Q_t^{zv}(x)}{Q_t^{zv}(u)} = \Delta_t^{zv}(x,u).
$$

Let now $\boldsymbol{\Delta}_1 = \{ \Delta_t^z, z \in X^{\Lambda \backslash \{t\}}, t \in \Lambda\}$ be a one-point transition energy field which elements satisfy conditions~\eqref{Delta-Markov}, and let $P_\Lambda$ be a corresponding to it random field with one-point conditional distribution $\boldsymbol{Q}_1(P_\Lambda)$. Then for all $t \in \Lambda$ and $x,u \in X^t$, $z \in X^{\partial t}$, $y,v \in X^{\Lambda \backslash (\{t\} \cup \partial t)}$, we have
$$
Q_t^{zy}(x) = \dfrac{\exp\{\Delta_t^{zy}(x,u)\}}{\sum \limits_{\alpha \in X^t} \exp\{\Delta_t^{zy}(\alpha,u)\}} = \dfrac{\exp\{\Delta_t^{zv}(x,u)\}}{\sum \limits_{\alpha \in X^t} \exp\{\Delta_t^{zv}(\alpha,u)\}} = Q_t^{zv}(x).
$$
Hence, $P_\Lambda$ is a Markov random field.
\end{proof}

Therefore, the following statement is true.

\begin{T}
\label{Th-Markov-Gibbs}
Any Markov random field $P_\Lambda$ is Gibbsian with the transition energy field $\boldsymbol{\Delta}_1$ satisfying conditions~\eqref{Delta-Markov}. Conversely, if the elements of the transition energy field $\boldsymbol{\Delta}_1$ of a (Gibbs) random field $P_\Lambda$ satisfy the conditions~\eqref{Delta-Markov}, then $P_\Lambda$ is a Markov random field.
\end{T}

It is not difficult to see that the well-known Hammersley--Clifford theorem on the equivalence of Markov finite random field and Gibbs finite random field with a pair potential of the nearest neighbors interactions is a particular case of Theorem~\ref{Th-Markov-Gibbs}. Indeed, according to Theorem~\ref{Th-allP-Gibbs}, Markov finite random fields are Gibbsian (in this case, there is no need to construct a potential). Conversely, let a neighborhood system $\partial = \{\partial t, t \in \Lambda\}$ in $\Lambda$ be given, and let $\Phi = \{\Phi_V, V \subset \Lambda\}$ be a pair potential of the nearest neighbors interactions, that is, $\Phi_{t \cup V} \neq 0$ only for $V \subseteq \partial t$, $t \in \Lambda$. For all $t \in \Lambda$ and $x,u \in X^t$, $z \in X^{\Lambda \backslash \{t\}}$, put
$$
\Delta_t^z(x,u) = \sum \limits_{\emptyset \neq J \subseteq \Lambda } \left(\Phi_{t \cup J}(u z_J) -  \Phi_{t \cup J}(x z_J) \right),
$$
and let $P_\Lambda^\Phi$ be a (Gibbs) random field corresponding to the transition energy field $\boldsymbol{\Delta}_1\Phi = \{\Delta_t^z, z \in X^{\Lambda \backslash \{t\}}, t \in \Lambda\}$. Since the conditions~\eqref{Delta-Markov} are satisfied, $P_\Lambda^\Phi$ is a Markov random field.

Let us also note that if  $H_1 = \{ H_t^z, z \in X^{\Lambda \backslash \{t\}}, t \in \Lambda\}$ is a Hamiltonian corresponding to a Markov random field $P_\Lambda$, then its elements satisfy the following conditions
$$
H_t^{zy}(u) - H_t^{zy}(x) = H_t^{zv}(u) - H_t^{zv}(x)
$$
for all $t \in \Lambda$ and $x,u \in X^t$, $z \in X^{\partial t}$, $y,v \in X^{\Lambda \backslash (t \cup \partial t)}$.

\section{Concluding remarks}

Thus, the methods developed for infinite random fields can be used to study random fields in finite volumes.

The main object considered in the paper is the system $\boldsymbol{\rm Q}_1$ of one-point probability distributions parameterized by boundary conditions and consistent in the sense~\eqref{Q1-sogl}. We showed that there is a one-to-one correspondence between the system $\boldsymbol{\rm Q}_1$ and compatible with it finite random field $P_\Lambda$. Further, we showed that there is a one-to-one correspondence between $\boldsymbol{\rm Q}_1$ and the one-point transition energy field $\boldsymbol{\Delta}_1$. These connections can schematically be presented by the following commutative diagram.

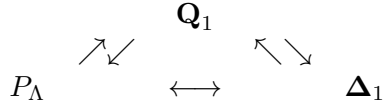
\begin{center}
\begin{tabular}{ccccccc}
   &  & $\boldsymbol{\rm Q}_1$ & & \\
 & $\nearrow \swarrow$ &   & $\nwarrow \searrow$  &   \\
 $P_\Lambda$ & & $\longleftrightarrow$ & & $\boldsymbol{\Delta}_1$
\end{tabular}
\end{center}

%We obtained (Theorem~\ref{Th-Q1->P+}) necessary and sufficient conditions for the existence and uniqueness of a finite random field $P_\Lambda$ with a given system of one-point conditional distributions $\boldsymbol{\rm Q}_1$. These conditions have a rather simple form and are expressed directly in terms of the elements of $\boldsymbol{\rm Q}_1$. Moreover, we suggest an (explicit) formula (see formula~\eqref{Q1->P}) connecting $P_\Lambda$ with the elements of $\boldsymbol{\rm Q}_1$ for all finite sets $\Lambda$.

One-to-one correspondence between $P_\Lambda$ and $\boldsymbol{\Delta}_1$ allows one to give the definition of a finite Gibbs random field, which does not involve the notion of potential.
From the condition~\eqref{Delta1Lambda-1} on the elements of $\boldsymbol{\Delta}_1$, it follows that each function $\Delta_t^z$ can be represented in the form
$$
\Delta_t^z(x,u) = H_t^z(u) - H_t^z(x), \qquad x,u \in X^t,
$$
where $H_t^z$ can be interpreted as a potential energy (Hamiltonian) of the system at point $t \in \Lambda$ under boundary condition $z \in X^{\Lambda \backslash \{t\}}$. The consistency conditions~\eqref{Delta1Lambda-2} impose certain conditions on the functions $H_t^z$, which lead to the presented axiomatic definition of a one-point Hamiltonian $H_1$ in a finite volume $\Lambda$. According to this definition, for any random field $P_\Lambda$, the set $H_1^\theta(P_\Lambda)$ of functions
$$
H_t^z(x) = \ln \frac{P_\Lambda(xz)}{P_\Lambda(\theta_t z)}, \qquad x \in X^t, z \in X^{\Lambda \backslash \{t\}}, t \in \Lambda,
$$
where $\theta \in X^\Lambda$ is some fixed configuration, is a Hamiltonian, and there is no need to represent $H_t^z$ as a sum of potentials.

The approach used made it possible to prove the fact that all finite random fields are Gibbsian.

Let us emphasize that obtained results do not depend on the structure of the finite set $\Lambda$. Its structure may be essential when considering special classes of finite random fields, for example, Markov finite random fields.

%\subsection{Further generalizations}

Our results carry over in a natural way to the case of infinite (both countable and continuous) measurable spaces $X$. Let, for example, $(X, \Im, \mu)$ be a complete separable metric space, where $X$ is a non-empty set, $\Im$ is a $\sigma$-algebra generated by open subsets of $X$, and $\mu$ is not identically zero countably additive measure.

%Further, let $(X^t, \Im^t, \mu_t)$, $t \in \Lambda$, be copies of the space $(X, \Im, \mu)$.

Consider a probability $P_\Lambda$ (random field) on the product space $(X^\Lambda, \Im^\Lambda, \mu_\Lambda)$ which has a positive density $p_\Lambda$ with respect to the measure $\mu_\Lambda$. For such random field, the system $\boldsymbol{Q}_1(P_\Lambda) = \{q_t^z, z \in X^{\Lambda \backslash \{t\}}, t \in \Lambda\}$ of its conditional densities is defined. In the case under consideration, the following result is an analogue of Theorems~\ref{Th-Q1->P+} and~\ref{Th-Q1-Delta1}.

\begin{T}
For a given system $\boldsymbol{\rm Q}_1 = \{q_t^z, z \in X^{\Lambda \backslash \{t\}}, t \in \Lambda \}$ of positive one-point probability densities parameterized by boundary conditions, there exists a unique probability density $p_\Lambda$ compatible with it if and only if the elements of $\boldsymbol{\rm Q}_1$ satisfy the consistency conditions~\eqref{Q1-sogl} and for any $u \in X^\Lambda$,
$$
Z_\Lambda(u) = \int \limits_{X^\Lambda} \prod \limits_{j=1}^n \frac{q_{t_j}^{(\alpha u)_j} (\alpha_j)}{q_{t_j}^{(\alpha u)_j} (u_j)} \mu_\Lambda(d\alpha) < \infty,
$$
where $\Lambda = \{t_1, t_2, ..., t_n\}$ is some enumeration of points in $\Lambda$, $n = \vert \Lambda \vert$. In this case, the elements of $\boldsymbol{\rm Q}_1$ have the Gibbs form with uniquely determined transition energy field $\boldsymbol{\Delta}_1^\Lambda$.
\end{T}

Note also, that the obtained results can be extended on the case when for the construction of the space $(X^\Lambda, \Im^\Lambda, \mu_\Lambda)$, to each point $t \in \Lambda$ corresponds a space $(X^t, \Im^t, \mu_t)$ which is not a copy of the same space $(X, \Im, \mu)$.

\textbf{Funding.} The work was supported by the Science Committee of RA in the frames of the research project 21AG-1A045.

\end{document}